\documentclass[11pt,a4paper]{article}
\usepackage{amsmath,amssymb,latexsym,amsthm}  
\usepackage[english]{babel}
\usepackage[latin2]{inputenc}

\begin{document}

\title{Some numerical characteristics of Sylvester and Hadamard matrices}
\author{\'Agota Figula and V. Kvaratskhelia} 
\date{}
\maketitle

\newtheorem{theor}{Theorem}[section]
\newtheorem{lem}[theor]{Lemma}
\newtheorem{co}[theor]{Corollary}
\newtheorem{prop}[theor]{Proposition}
\newtheorem{remark}[theor]{Remark}
\newtheorem{conj}[theor]{Conjecture}
\newtheorem{example}[theor]{Example}

\newcommand{\nc}{\newcommand}
\newcommand{\Pro}{\noindent{\em Proof.}}
\renewcommand{\th}{\vartheta}
\nc{\li}{\limits} \nc{\ri}{\right} \nc{\ph}{\varphi}
\nc{\la}{\lambda} \nc{\al}{\alpha} \nc{\eps}{\varepsilon}
\nc{\e}{\mathbb E} \nc{\del}{\Delta} \nc{\Om}{\Omega}
\nc{\om}{\omega} \nc{\gm}{\gamma} \nc{\ld}{\ldots} \nc{\cd}{\cdots}
\nc{\de}{\delta} \nc{\ro}{\varrho} \nc{\be}{\beta}
\nc{\Gm}{\varGamma} \nc{\sg}{\sigma} \nc{\ka}{\kappa}
\nc{\ov}{\overline} \nc{\lr}{\Leftrightarrow} \nc{\ra}{\Rightarrow}
\nc{\lag}{\langle} \nc{\rag}{\rangle} \nc{\lf}{\left}
\nc{\sgn}{\operatorname{sgn}}

\begin{abstract} We introduce numerical characteristics of
{\it Sylvester} and {\it Hadamard matrices} and
give their estimates and some of their applications.
\end{abstract}


\noindent
{\footnotesize {2010 {\em Keywords:} Hadamard matrices;
Sylvester matrices; Banach space; subsymmetric basis.}}

\noindent
{\footnotesize {\em Acknowledgments:}  This paper are supported by the European Union's Seventh
Framework Programme (FP7/2007-2013) under grant agreements no. 317721, no. 318202,
by the Shota Rustaveli National Science Foundation
grant no. FR/539/5-100/13 and by the J\'anos Bolyai Research Fellowship.}

\section{Introduction}

\noindent
A Hadamard matrix has a simple structure, it is a square matrix such that its entries are either $+1$ or $-1$ and its rows (columns)
are mutually orthogonal. In spite of the fact that Hadamard matrices have been actively studied for
about 150 years, they still have unknown properties. If $\mathcal H_n$ is a Hadamard matrix of order $n$, then the matrix
$$\left[\begin{matrix} \mathcal H_n & \mathcal H_n\\
\mathcal H_n &- \mathcal H_n\end{matrix}\right]$$
is a Hadamard matrix of order $2n$. Applying this algorithm repeatedly J.J. Sylvester has constructed a particular sequence of
Hadamard matrices of order $2^n$. These matrices are called Sylvester matrices or Walsh matrices.

Hadamard matrices have a wide range of applications in the code theory, scheduling theory, statistics,
modern communications etc. In this paper we deal with its application for certain problems of
functional analysis. Namely, using Hadamard matrices, in the classical
Banach spaces it is easy to construct examples of unconditionally convergent series which do not converge absolutely
(see \cite{BS}, \cite{VK1}, \cite{VK2}). Note that in \cite{MP} the author made a considerable
effort to prove the existence of such series in the space $l_1$ without giving construction.
To prove the unconditional convergence of the above mentioned constructed series the numerical characteristics of 
Hadamard and Sylvester matrices are important tools. In the present paper the general forms of
these tools for Banach spaces with bases  are considered. These
characteristics and the structure of the Hadamard and
Sylvester matrices play an important role in the investigation of the convergence of series in Banach
spaces (see e.g., \cite{C-T}, \cite{VK3}, \cite{VK4}). For these characteristics we give estimates
(cf. Theorems 3.1, 3.6, 4.2 and 4.8). We believe that the investigated characteristics and their estimates complete 
our knowledge about Hadamard matrices and may have applications in other fields of mathematics.

In Section 2 some concepts, definitions and auxiliary results required for further discussions are given.

In Section 3 the numerical characteristic $\ro^{(n)}$ of Sylvester matrices is introduced
and its estimates for the case of a Banach space with a subsymmetric basis $(\varphi_i)$ are studied.
For every positive integer $n$ we prove the following estimates (cf. Theorems 3.1 and 3.6)
$$\max\left\{{\frac{n+2}{6}} \cdot \lambda(2^n), \, 2^n\right\} \le \ro^{(n)}
\le \min\left\{\left(1+\sum\limits_{j=1}^n2^{-j}\lambda(2^{j-1})
\right)\cdot 2^n, \lambda(n)\cdot 2^n\right\},$$
where $\la(n)=\Vert\sum\li_{i=1}^n\ph_i\Vert$.

In Section 4 we define the analogue characteristic $\varrho_n$ for Hadamard
matrices. For every positive integer $n$ for which there is a Hadamard matrix
of order $n$ we show the following estimates (cf. Theorems 4.2 and 4.8)
$$\max\left\{(1/\sqrt2)\,\lambda(n)\,\sqrt n, \, n\right\}
\le \varrho_n \le \lambda([\sqrt n]+1)\,n, $$
where $[\sqrt n]$ is the integer part of $\sqrt n$.

As an application of the introduced notions we give a characterization for the spaces isomorphic
to $l_1$ in terms of these characteristics (cf. Theorem 4.10).

In Section 5 we pose an open problem which has naturally arisen from our investigations.

Most of the results of this paper were announced in \cite{Kv1} without proofs.
Here these results and some new ones are given with complete proofs.

\section{Notation and Preliminaries}

We follow the standard notation and terminology used, for example, in \cite{LT}.
The notations $c_0, l_p$ and $L_p, 1\le p<\infty,$ have their usual meaning.

A sequence $(\ph_i)$ of nonzero elements in a real Banach space $X$ is called a ({\it Schauder})
{\it basis} of $X$ if for every $x\in X$ there is a unique sequence of
scalars $(\al_i)$ so that $x=\sum\li_{i=1}^\infty\al_i\ph_i$. If $(\ph_i)$
is a basis in a Banach space $X$ with a norm $\|\cdot\|$, then there is
a constant $K\ge 1$ so that for every choice of scalars $(\al_i)$
and positive integers $n<m$, we have
$$\|\sum\li_{i=1}^n\al_i\ph_i\|\le K \|\sum\li_{i=1}^m\al_i\ph_i\|.$$
The smallest possible constant $K$ in this inequality is called the {\it basis
constant} of $(\ph_i)$. Note that in $X$ there exists an
{\it equivalent norm} $|||\cdot|||$ (i.e. for some positive constants
$C_1, C_2$: $C_1\|x\|\le|||x|||\le C_2\|x\|$ for every $x\in X$)
for which the basis constant is $K=1$.

A basis $(\ph_i)$ is called {\it normalized} if $\|\ph_i\|=1$ for all $i$.
Let $(\ph_i)$ be a basis of a Banach space $X$.
A sequence of linear bounded functionals $(\ph_i^*)$ defined by the relation
$\langle \ph_i^*,\ph_j\rangle=\delta_{ij}$, where $\delta_{ij}$ is the
{\it Kronecker delta}, is called the sequence of {\it biorthogonal functionals}
associated to the basis $(\ph_i)$. Two bases, $(\ph_i)$ of $X$ and $(\psi_i)$
of $Y$, are called equivalent provided a series $\sum\li_{i=1}^\infty\al_i\ph_i$
converges if and only if $\sum\li_{i=1}^\infty\al_i\psi_i$ converges.

A basis $(\ph_i)$ of a Banach space $X$ is {\it unconditional}
if for any permutation $\pi:\mathbb{N}\to \mathbb{N}$ of the set $\mathbb N$ of positive integers
$(\ph_{\pi(i)})$ is a basis of $X$. If $(\ph_i)$ is an unconditional basis of a real Banach space $X$, then there is a
constant $K\ge 1$ so that for every choice of scalars $(\al_i)$ for which $\sum\li_{i=1}^\infty\al_i\ph_i$
converges and every choice of bounded scalars $(\la_i)$ we have
$$\|\sum\li_{i=1}^\infty\la_i\al_i\ph_i\|\le K \sup\li_i |\la_i| \|\sum\li_{i=1}^\infty\al_i\ph_i\|.$$
The smallest possible constant $K$ in this inequality is called the {\it unconditional
constant} of $(\ph_i)$. If $(\ph_i)$ is an unconditional basis of $X$, then
there is an equivalent norm in $X$ so that the unconditional constant becomes 1.

The sequence of unit vectors $e_i=(0,0,\ld,\buildrel{i}\over{1},0,\ld), i=1,2,\ld,$
is an example of an unconditional basis in $c_0$ and $l_p,\,1\le p<\infty$
(the basis $(e_i)$ is called the {\it natural basis} of the corresponding spaces).
The {\it Haar system} is an unconditional basis in the function spaces $L_p(0,1),\,1< p<\infty$.
This system is also basis in $L_1(0,1)$, but in this space there does not exist an
unconditional basis.

Every normalized unconditional basis in $l_1, l_2$ or $c_0$ is equivalent to the
natural basis of these spaces. Moreover, a Banach space has, up to equivalence, a unique
unconditional basis if and only if it is isomorphic to one of the following three
spaces: $l_1, l_2$ or $c_0$.

Let $(X,\|\cdot\|)$ be a Banach space with a normalized basis $(\ph_i)$. Consider the expression
$$\la(n)=\Vert\sum\li_{i=1}^n\ph_i\Vert, \quad n=1,2,\ldots .$$
For every space having an unconditional basis whose unconditional constant is $1$ with the exception of the space $c_0$
we have that $(\la(n))$ is a non-decreasing sequence and $\lim\li_{n\to\infty}\la(n)=\infty$. More precisely, if
$\sup\li_n\la(n)<\infty$, then $(\ph_i)$
is equivalent to the natural basis of the space $c_0$ (see, for example, \cite{LT}, p. 120).

A basis $(\ph_i)$ of a Banach space $X$ is said to be {\it symmetric} if for any permutation
$\pi$ of the positive integers $(\ph_{\pi(i)})$ is equivalent to $(\ph_i)$. If
$(\ph_i)$ is a symmetric basis of a Banach space $X$, then there is a
constant $K$ such that for any choice of scalars $(\al_i)$ for which $\sum\li_{i=1}^\infty\al_i\ph_i$
converges, every choice of signs $\vartheta=(\vartheta_i)$ and any permutation
$\pi$ of the positive integers we have
$$\|\sum\li_{i=1}^\infty\vartheta_i\al_i\ph_{\pi(i)}\|\le K \|\sum\li_{i=1}^\infty\al_i\ph_i\|.$$
The smallest possible constant $K$ in this inequality is called the {\it symmetric
constant} of $(\ph_i)$.

A basis $(\ph_i)$ of a Banach space $X$ is called {\it subsymmetric} if it is unconditional
and for every increasing sequence of integers $(i_n)$, $(\ph_{i_n})$ is equivalent to $(\ph_i)$.
If $(\ph_i)$ is a subsymmetric basis of a Banach space $X$, then there is a
constant $K$ such that for any choice of scalars $(\al_i)$ for which $\sum\li_{i=1}^\infty\al_i\ph_i$
converges, every choice of signs $\vartheta=(\vartheta_i)$ and every
increasing sequence of integers $(i_n)$ we have
$$\|\sum\li_{n=1}^\infty\vartheta_n\al_n\ph_{i_n}\|\le K \|\sum\li_{i=1}^\infty\al_i\ph_i\|.$$
The smallest possible constant $K$ in this inequality is called the {\it subsymmetric
constant} of $(\ph_i)$.

Every symmetric basis is subsymmetric. The converse of this assertion is not true.
The unit vectors in $l_p,\, 1\le p<\infty,$ and $c_0$ are examples of symmetric basis.
\vskip 0.2cm

\begin{prop}\label{2.1}
(see \cite{LT}, Proposition 3.a.7, p. 119). Let $(X,\|\cdot\|)$
be a Banach space with a symmetric basis $(\ph_i)$ whose symmetric
constant is  equal to 1. Then there exists a new norm $\|\cdot\|_0$ on $X$ such that:

$(a).$ $\|x\|\le\|x\|_0\le 2\|x\|$ for all $x\in X$;

$(b).$ The symmetric constant of $(\ph_i)$ with
respect to $\|\cdot\|_0$ is equal to 1;

$(c).$ If we put $\la_0(n)=\Vert\sum\li_{i=1}^n\ph_i\Vert_0,
n=1,2,\ldots ,$ then $\{\la_0(n+1)-\la_0(n)\}$ is a non-increasing sequence, i.e. $\la_0(\cdot)$
is a concave function on the integers.
\end{prop}
\vskip 0.2cm

The converse of the last assertion is also true
in the sense that, for every concave non-decreasing
sequence of positive numbers $(\la_k)$ there exists at least one Banach space $X$ having a
symmetric basis $(\ph_i)$ with symmetric constant equal to 1 such that
$\Vert\sum\li_{i=1}^n\ph_i\Vert=\la_n$ for every $n$.
\vskip 0.2cm

\begin{prop}\label{2.2}
(see \cite{LT}, Proposition 3.a.4, p. 116).  ${\bf (A).}$ Let $X$ be a Banach space
with a normalized subsymmetric basis $(\ph_i)$ whose subsymmetric
constant is 1. Then the following inequality is valid
$$
\Vert\sum\li_{i=1}^n\al_i\ph_i\Vert\ge {\frac{\sum\li_{i=1}^n|\al_i|}{n}}
\, \la(n), \quad n=1,2,\ldots .
$$

${\bf (B).}$ Moreover, if $(\ph_i)$ is a subsymmetric basis, then one has

$$
\Vert\sum\li_{i=1}^n\al_i\ph_i\Vert\ge {\frac{\sum\li_{i=1}^n|\al_i|}{2n}}
\, \la(n), \quad n=1,2,\ldots .
$$
\end{prop}
\vskip 0.2cm

From this it follows that if $\lim\li_{n\to\infty}\sup\la(n)/n>0$, then
$(\ph_i)$ is equivalent to the natural basis of the space $l_1$
(see, for example, \cite{LT}, p. 120).

The {\it Rademacher functions} $r_k, \, k=1,2,\ldots$, are defined on $[0,1]$ by the equality
$$r_k(t)=sign(\sin{2^k{\pi t}}).$$

Let us note the well-known {\it Khintchine's inequality}: for every $0<p<\infty$
there exist positive constants $A_p$ and $B_p$ so that
$$A_p\lf(\sum\li_{k=1}^m |\al_k|^2\ri)^{1/2}\le \lf(\int\li_0^1
\lf|\sum\li_{k=1}^m\al_k r_k(t)\ri|^p\,dt\ri)
^{1/p}\le B_p\lf(\sum\li_{k=1}^m |\al_k|^2\ri)^{1/2}, $$
$m=1,2,\ldots,$ for every choice of scalars $(\al_1,\al_2,\ldots,\al_m)$.
For $p=1$ the best constant is $A_1=1/\sqrt2$ (see \cite{Sz}).

A Banach space $X$ is said to be of {\it type $p$} if there is a constant
$T_p=T_p(X)\ge 0$ such that for any finite collection of vectors
$x_1, x_2, \ld, x_n$ in $X$ we have
$$\lf(\int_0^1\lf\Vert\sum\li_{k=1}^n r_k(t) x_k\ri\Vert^2\,dt\ri)^{1/2}\le T_p\,
\lf(\sum\li_{k=1}^n\lf\Vert x_k\ri\Vert^p\ri)^{1/p}, \quad n=1,2,\ldots .$$

In the Khintchine's inequality the notion of type $p$ has meaning for the case
$0<p\le 2$. Every Banach space has type $p$ for $0<p\le 1$.
The spaces $l_p, L_p([0,1]), 1\le p<\infty$, have type
$\min(2,p)$.

{\it A Hadamard matrix} is a square matrix of order $n$ with entries
$\pm1$ such that any two columns (rows) are orthogonal (see {\it
e.g.} \cite{Hl}, p. 238, \cite{WS}, p. 44). We denote by $\mathcal H_n=\left[h_{ki}^n\right]$
a Hadamard matrix of order $n$.  It is
easy to see that the order of a Hadamard matrix is either $1$ or $2$ or it is
divisible by 4. {\it Hadamard} put forward the conjecture that
for any $n$ divisible by $4$ there exists a Hadamard matrix of
order $n$. As far as we know, {\it Hadamard's conjecture} remains
open. Let $\mathbb{N}_{\mathcal H}$ be the set of all positive integers
$n$ for which there exists a Hadamard matrix of order $n$.

The following property follows from the definition of Hadamard matrices.
If $\mathcal H_n=\left[h_{ki}^n\right]$ is a Hadamard matrix, then for every
$n,\,n\in\mathbb{N}_{\mathcal H}$, we have
$$\sum\li_{i=1}^n h_{ki}^n h_{mi}^n=n\,\de_{km},\quad \sum\li_{k=1}^n h_{ki}^n h_{kj}^n=n\,\de_{ij}.$$
Therefore for any $n,\,n\in\mathbb{N}_{\mathcal H}$, and every sequence $(\be_i)_{i\le n}$ of
real numbers one has
$$\sum\li_{k=1}^n\lf(\sum\li_{i=1}^n h_{ki}^n\be_i\ri)^2=n\sum\li_{i=1}^n\be_i^2.$$

It is easy to see that multiplying any row or any column of a Hadamard matrix by $-1$
we get again a Hadamard matrix.

Let the triple $(\Om,\mathfrak A,\mathbb{P})$ be a probability space,
where $\Om$ be a non-empty set, $\mathfrak A$ be a $\sigma$-algebra of
subsets of $\Om$ and $\mathbb{P}$ be a probability measure on the measurable space
$(\Om,\mathfrak A)$, (i.e. $\mathbb{P}$ is assumed to be a non-negative measure on
$(\Om,\mathfrak A)$ satisfying the condition $\mathbb{P}(\Om)=1$). Let $X$ be
a real Banach space with the topological dual space $X^*$. A function
$\xi:\Om\to X$ is scalarly measurable (respectively scalarly integrable)
if for each $x^*\in X^*$ the scalar function $\langle x^*,\xi\rangle$ is measurable
(respectively integrable, i.e. $\langle x^*,\xi\rangle\in L_1(\Om,\mathfrak A,\mathbb{P})$).
A scalarly integrable function $\xi:\Om\to X$ is {\it Pettis integrable}
(or {\it weak integrable}) if for each $A\in\mathfrak A$ there exists a vector $m_{\xi,A}\in X$ such that for every $x^*\in X^*$ we have
$$\langle x^*,m_{\xi,A} \rangle=\int\li_A \langle x^*,\xi \rangle \,d \mathbb{P}.$$
For a Pettis integrable function $\xi:\Om\to X$ the element $m_{\xi,\Om}$
is called the {\it Pettis integral} of $\xi$ with respect to $\mathbb{P}$. It is also called
the {\it mean value} of the function $\xi$. We denote by $\mathbb{E}\,\xi$ the Pettis integral of the function $\xi$. If a function
$\xi:\Om\to X$ has a measurable norm and
there exists $\mathbb{E}\,\xi$, then $\Vert\mathbb{E}\,\xi\Vert\le\mathbb{E}\,\Vert\xi\Vert$.
For every separably valued function $\xi:\Om\to X$ from the condition
$\mathbb{E}\,\Vert\xi\Vert<\infty$ it follows the existence of the Pettis integral
$\mathbb{E}\,\xi$ ($\xi$ is separably valued if $\xi(\Om)$ is a separable subset of $X$).

For details and proofs related with the topics of this section see \cite{LT} and \cite{VTC}.

\section{Sylvester matrices}

The {\it Sylvester matrices} are special cases of Hadamard matrices.
They are defined by the recursion relations (cf. \cite{WS}, p. 45):
$${\mathcal S}^{(1)}=\left[\begin{matrix}1&1\\1&{-1}\end{matrix}\right],
\quad
{\mathcal S}^{(n)}=\left[\begin{matrix}{\mathcal S}^{(n-1)}&{\mathcal S}^{(n-1)}\\
{\mathcal S}^{(n-1)}&-{\mathcal
S}^{(n-1)}\end{matrix}\right],\quad n=2,3,\ldots.$$
${\mathcal S}^{(n)}$ is a Hadamard matrix of order
$2^n$ and hence $2^n\in\mathbb{N}_{\mathcal H}$ for all $n=1,2, \ldots$.

If the first column of a Hadamard matrix $\mathcal
H_n=\left[h_{ki}^n\right]$ consists of only $+1$, then one has
$$\sum\li_{k=1}^n h_{ki}^n=
\begin{cases} \ \ n,& \ \ \ {\hbox{for}} \ \ \ i=1, \cr
\ \ 0,& \ \ \ {\hbox{for}} \ \ \ i=2,3,\ldots ,n. \cr
\end{cases}$$

In particular, if $\mathcal S^{(n)}=\left[s_{ki}^{(n)}\right]$ is the Sylvester matrix of order $2^n$,
$n=1,2,\ldots$, then we get
$$\sum\li_{k=1}^{2^n}s_{ki}^{(n)}=
\begin{cases} \ \ 2^n,& \ \ \ {\hbox{for}} \ \ \ i=1, \cr
\ \ 0,& \ \ \ {\hbox{for}} \ \ \ i=2,3,\ldots ,2^n \cr
\end{cases}$$
and
{\hsize 17.25cm {$$\sum\li_{k=1}^{2^{n-1}}s_{ki}^{(n)}=
\begin{cases} \ \ 2^{n-1},& \ \ \ {\hbox{for}} \ \ \ i=1 \ \ \
{\hbox{and}} \ \ \ i=2^{n-1}+1, \cr \ \ 0,& \ \ \ {\hbox{otherwise.}} \cr
\end{cases}$$}}

Let ${\mathcal S}^{(n)}=\left[s_{ki}^{(n)}\right]$ be the Sylvester matrix of order $2^n$,
$n=1,2,\ldots$, and $X$ be a Banach space
with a norm $\|\cdot\|$ and a normalized basis $(\varphi_i)$.
Consider the function
$$\varrho^{(n)}(m)=\left\Vert\sum\limits_{i=1}^{2^n}\left(\sum\limits_{k=1}^m
s_{ki}^{(n)}\right) \varphi_i\right\Vert,\, m=1,2,\ldots,2^n. \eqno (3.1) $$
One has $\ro^{(n)}(1)=\lambda(2^n),
\ro^{(n)}(2)=2\,\lf\Vert\sum\li_{i=1}^{2^{n-1}}\ph_{2i-1}\ri\Vert,
\ro^{(n)}(2^n)=2^n,$ where
$\lambda(2^n)=\lf\Vert\sum\li_{i=1}^{2^n}\ph_i\ri\Vert$. The
function $\varrho^{(n)}(m)$ obviously depends on $X$, the norm
in $X$ and the choice of basis $(\varphi_i)$. In particular,
for the case of the spaces $l_p, \, 1\le p<\infty$, with respect to the
natural basis $\ro^{(n)}(m)$ has the form
$\left(\sum\li_{i=1}^{2^n}\lf|\sum\li_{k=1}^m
s_{ki}^{(n)}\ri|^p \ri)^{1/p}.$

We set
$$\varrho^{(n)}=\max\limits_{1\le m\le 2^n}\varrho^{(n)}(m). \eqno (3.2)$$

The function $\varrho^{(n)}(m)$ can be expressed  as follows.
Let $a_k=\sum\limits_{i=1}^{2^n} s_{ki}^{(n)}\varphi_i, \
k=1,2,\ldots,2^n.$ Then one has
$\varrho^{(n)}(m)=\left\Vert\sum\limits_{k=1}^m a_k\right\Vert.$
If $(\varphi_i)$ is an unconditional basis with unconditional
constant equal to $1$, then, obviously, $||a_k||=\lambda(2^n)$ for any
$k=1,2,\ldots,2^n$ and $\varrho^{(n)}\le\lambda(2^n) \,
2^n\le2^{2n}$.

In $l_p, \ 1\le p<\infty,$ it was proved in \cite{VK1} that
$\ro^{(n)}\le n\,2^n.$

The following theorem gives a similar estimate of $\varrho^{(n)}$ in the
case of general Banach spaces with subsymmetric basis.
\vskip 0.2cm

\begin{theor}\label{3.1}
Let $X$ be a Banach space with normalized subsymmetric basis whose
subsymmetric constant is $1$. Then for $\varrho^{(n)}$ defined by
(3.2) one has the following estimate
$$\varrho^{(n)}\le \min\left\{\left(1+\sum\limits_{j=1}^n2^{-j}\lambda(2^{j-1})
\right)\cdot 2^n, \quad \lambda(n)\cdot 2^n\right\} , \quad
n=1,2,\ldots. \eqno (3.3)$$
\end{theor}
\vskip 0.2cm

\Pro \
First we prove the inequality $\varrho^{(n)}\le \left(1+\sum\limits_{j=1}^n2^{-j}
\lambda(2^{j-1})\right)\cdot 2^n$ by induction.
For $n=1$ it is true since the left hand side of (3.3) is equal to $2$ and the right
hand side is equal to $3$. Let $n\ge 2$. Introduce the following notation
$$\al_i^{(n)}(m)=\sum\li_{k=1}^m s_{ki}^{(n)}, \quad 1\le i, m\le 2^n.
\eqno (3.4)$$

Therefore we get
$$\al_1^{(n)}(m)=m \eqno(3.5)$$
and
$$\al_{2^{n-1}+1}^{(n)}(m)=\begin{cases} m,& \ {\hbox{for}} \quad
1\le m\le 2^{n-1},\cr 2^n-m,& \ {\hbox{for}} \quad 2^{n-1}+1\le
m\le 2^n.
\end{cases} \eqno(3.6)$$
Since $i\le2^n$ we can write that
$i=\eps_n2^n+\eps_{n-1}2^{n-1}+\ldots + \eps_12+\eps_0,$ where
$\eps_j \in \{0,1\}$ for every $j$. Then by the definition
and the properties of the Sylvester matrices we can
prove by induction that for any $i$
$$ \max\li_{1\le m\le 2^n} \lf|\al_i^{(n)}(m)\ri|=2^{f(i)}, \eqno (3.7)$$
where the function $f:\{1,2,\ldots,n\}\to \{0,1,2,\ldots,n\}$ is defined as follows: 
$f(1)=n$; $f(i)=0$ if $\eps_0=0$ (i.e. $i$ is an even number)
and if $\eps_0=1$ (i.e. $i$ is an odd number), then for $f(i)$ we have:
$\eps_{f(i)}=1$ and $\eps_j=0$ for every $j=1,2,\ldots,f(i)-1.$

For $i=1$ and $i=2^{n-1}+1$ the equality (3.7) is valid
since from the relations (3.5) and (3.6) it follows that
$\max\li_{1\le m\le2^n}\lf|\al_1^{(n)}(m)\ri|=2^n$ and
$\max\li_{1\le m\le2^n}\lf|\al_{2^{n-1}+1}^{(n)}(m)\ri|=2^{n-1}$.
To prove (3.7) for the rest indexes $i$ we use the following equalities
$$\max\li_{1\le m\le2^{n+1}}\lf|\al_{2^n+i}^{(n+1)}(m)\ri|=
\max\li_{1\le m\le2^{n+1}} \lf|\al_i^{(n+1)}(m)\ri|=\max\li_{1\le
m\le2^n}\lf|\al_i^{(n)}(m)\ri| \eqno (3.8)$$ for any
$i=2,3,\ldots,2^n$, which is a consequence of the definition
and the properties of the Sylvester matrices.
Every positive integer $i, \,1\le i\le2^{n+1},$ has the unique
representation given by
$$i=\begin{cases} \ \eps_n2^n+\ldots+\eps_12+\eps_0,& \ \ \ {\text{for}}
\ \ \ 1\le i\le2^n,\cr \ 2^n+\eps_n 2^n+\ldots+\eps_1 2+\eps_0,& \
\ \ {\text{for}} \ \ \ 2^n+1\le i\le2^{n+1}.\cr
\end{cases} \eqno (3.9)$$
If $i$ is an even number, then in (3.9) we have $\eps_0=0$
and by (3.7) and (3.8) we obtain $\max\li_{1\le
m\le2^{n+1}}\lf|\al_i^{(n+1)}(m)\ri|=1$. If $i$ is an odd number
and, in addition, $i\ne1$ and $i\ne 2^n+1$, then we can rewrite (3.9) as follows:
$$i=\begin{cases} \
\eps_n2^n+\ldots+\eps_{j_0+1}2^{j_0+1}+2^{j_0}+1,& \ \ \
{\text{for}} \ \ \ \ \ 3\le i\le2^n,\cr
2^n+\eps_n2^n+\ld+\eps_{j_0+1}2^{j_0+1}+2^{j_0}+1,& \ \ \
{\text{for}} \ \ \ \ \ 2^n+3\le i\le2^{n+1},\cr
\end{cases}$$
where $j_0=1,2,\ldots,n-1.$ Using again relations (3.7) and (3.8) we certainly have
$\max\li_{1\le
m\le2^{n+1}}\lf|\al_i^{(n+1)}(m) \ri|=2^{j_0}.$

Applying now a simple combinatorial calculation we get that the number
of indexes $i, \ 1\le i\le 2^n$, for which $\max\li_{1\le m\le 2^n} \lf|\al_i^{(n)}(m)\ri|=2^j$,
is equal to $2^{n-j-1}$ for $j=0,1,2,\ldots,n-1,$ and the equality
$\max\li_{1\le m\le 2^n}\lf|\al_i^{(n)}(m)\ri|=2^n$ is satisfied only for $i=1$.

As the subsymmetric constant of the basis $(\ph_i)$ is $1$, using (3.7), we obtain
for every $m=1,2,\ld,2^n$ the following relations:
$$\ro^{(n)}(m)=\lf\Vert\sum\li_{i=1}^{2^n}\lf|\al_i^{(n)}(m)\ri|\ph_i\ri\Vert\le
\lf\Vert\sum\li_{i=1}^{2^n}\max\li_{1\le m\le
2^n}\lf|\al_i^{(n)}(m)\ri|\ph_i\ri\Vert=$$
$$=\lf\Vert 2^n\ph_1+\sum\li_{j=1}^n2^{n-j}\sum\li_{i=2^{j-1}}^{2^j-1}\ph_
{\pi(i+1)}\ri\Vert, \eqno (3.10)$$ where $\pi$ is a permutation of a sequence of the positive integers
$\{2,3,\ldots,2^n\}$. Applying now the triangular inequality on the right hand side of
(3.10) and using the fact that $(\ph_i)$ is a subsymmetric basis we get the required inequality.

Now we prove the inequality $\varrho^{(n)}\le\lambda(n)\,2^n.$
The number of the (not necessarily different) basis elements involved in the right hand side
of the inequality (3.10) is equal or less than $n\cdot2^n$ (more exactly, $(1+n/2)\cdot2^n$).
Hence we get the following equality
$$2^n\ph_1+\sum\li_{j=1}^n2^{n-j}\sum\li_{i=2^{j-1}}^{2^j-1}\ph_{\pi(i+1)}=
\sum\li_{k=1}^{2^n}\sum\li_{i=1}^{l_k}\ph_{k_i}, \eqno (3.11)$$ where $1\le l_k\le n$
for any $k=1,2,\ldots, 2^n$, $\ph_{k_i}\in\left\{\ph_1,\ph_2,\ldots,\ph_{2^n}\right\}$,
for any fixed $k$ and for every $i\ne j, \, i,j=1,2,\ldots,l_k$, we have $\ph_{k_i}\ne\ph_{k_j}$
and for any fixed $i$ but for different indexes $k$ the elements $\ph_{k_i}$ can be the same.
As the basis $(\ph_i)$ is subsymmetric with subsymmetric constant equal to $1$
and (3.10) and (3.11) are valid we obtain
$$\ro^{(n)}(m)\le\lf\Vert2^n\ph_1+\sum\li_{j=1}^n2^{n-j}\sum\li_{i=2^{j-1}}^
{2^j-1}\ph_{\pi(i+1)}\ri\Vert=\lf\Vert\sum\li_{k=1}^{2^n}\sum\li_{i=1}^{l_k}
\ph_{k_i}\ri\Vert\le$$
$$\le\sum\li_{k=1}^{2^n}\lf\Vert\sum\li_{i=1}^{l_k}\ph_{k_i}\ri\Vert
=\sum\li_{k=1}^{2^n}\lambda(l_k)\le\lambda(n)\,2^n$$
for every $m$. This proves the theorem.
\qed
\vskip 0.2cm

\begin{remark}\label{3.2}
{\rm For the estimates proved
in Theorem \ref{3.1} with respect to the natural basis we obtain
the relation $1+\sum\limits_{j=1}^n2^{-j}\lambda(2^{j-1}) \le \lambda(n)$ in the case of $X=l_1$,
but we have  the converse relation
$1+\sum\limits_{j=1}^n2^{-j}\lambda(2^{j-1}) \ge \lambda(n)$ in the case of $X=c_0$.}
\end{remark}
\vskip 0.2cm

Let $X$ be a Banach space (not necessarily with basis),
$x_1,x_2,\ldots,x_{2^n}$ be a sequence of elements from the unit
ball of $X$ and ${\mathcal S}^{(n)}$ be the Sylvester matrix of
order $2^n,\,n=1,2,\ldots$. By analogy with the definition of
$\varrho^{(n)}$ let
$\hat\varrho^{(n)}(m)=\left\Vert\sum\limits_{i=1}^{2^n}\left(\sum\limits_{k=1}^ms_{ki}^{(n)}
\right)x_i\right\Vert, \, m=1,2,\ldots,2^n$, and let
$\hat\varrho^{(n)}=\max\limits_{1\le m\le
2^n}\hat\varrho^{(n)}(m).$
\vskip 0.2cm

\begin{co}\label{3.3}
 We have $\hat\varrho^{(n)}\le n\cdot2^n.$
 \end{co}

\Pro \
Using the triangular inequality and the fact that $||x_i||\le 1$ for any $i$, we have
$$\hat\ro^{(n)}(m)\le\sum\li_{i=1}^{2^n}\lf|\sum\li_{k=1}^ms_{ki}^{(n)}\ri|.$$
The right hand side of the last relation is the expression $\ro^{(n)}(m)$ in the space
$l_1$ with respect to the natural basis, which is for every $m=1,2,\ld,2^n$
less or equal than $n\,2^n$ (cf. Theorem \ref{3.1}).
\qed
\vskip 0.2cm

\begin{co}\label{3.4}
 Let $X$ be a Banach space of type $p, \, p>1,$ with a normalized
subsymmetric basis $(\varphi_i)$ whose subsymmetric constant is 1. Then one has
$$\varrho^{(n)}\le c \cdot 2^n,$$
where the constant $c\ge1$ depends only on the space $X$.
\end{co}

\Pro \
Since $(\ph_i)$ is a normalized subsymmetric basis whose subsymmetric constant is 1,
then $\lambda(2^{j-1})\le T_p(X) \,2^{(j-1)/p}$ for every $j\ge1$, where $T_p(X)$
is the constant involved in the definition of the space of type $p$. Then
for the right hand side of (3.3) we get
$$1+\sum\li_{j=1}^n2^{-j}\lambda(2^{j-1})\le{1+T_p(X)\sum\li_{j=1}^n2^{-j+(j-1)/p}}\le1+T_p(X)/(2-2^{1/p}).$$
Taking $c=1+T_p(X)/(2-2^{1/p})$ the proof is finished.
\qed
\vskip 0.2cm

Let us note that in the space $c_0$ we have a similar estimate,
namely $\ro^{(n)}\le 2^n$ (cf. Theorem \ref{3.1}), although $c_0$ is a space of type
1. As $\ro^{(n)}\ge2^n,$ we get $\ro^{(n)}=2^n$ in the space $c_0$.

Thus, in the Banach spaces of type $p,\,p>1,$ (as well as in $c_0$), we have
$\sup\li_n\varrho^{(n)}/2^n<\infty$. But in general this is not true.
The following statement shows the validity of this fact for the space $l_1$.
\vskip 0.2cm

\begin{theor}\label{3.5}
\cite{Kv2}. For the space $l_1$ with the natural basis one has
$$\varrho^{(n)}=\max\limits_{1\le m\le2^n}\varrho^{(n)}(m)={(3n+7)2^n}/9+{2(-1)^n}/9,
\quad n\ge 1.$$ For any $n$ the maximum is attained at
the points $m_n={(2^{n+1}+(-1)^n)}/3$ and $m_n^{'}={(5\cdot2^{n-1}+(-1)^{n-1})}/3$ .
\end{theor}
\vskip 0.2cm

Let us estimate $\ro^{(n)}$ from below.
\vskip 0.2cm

\begin{theor}\label{3.6}
If a Banach space $X$ satisfies the conditions of
Theorem \ref{3.1}, then one has
$$\ro^{(n)}\ge\max\left\{{\frac{n+2}{6}} \, \lambda(2^n), \ 2^n\right\},
\quad n=1,2,\ldots.$$
\end{theor}

\Pro \
By the definition of $\ro^{(n)}(m)$ for any positive integer $n$ we have,
$\ro^{(n)}\ge\ro^{(n)}(2^n)=2^n$ and the inequality $\ro^{(n)}\ge2^n$ is evident.

Let us prove that for any integer $n$ the inequality
$\ro^{(n)}\ge{\frac{n+2}{6}} \,\lambda(2^n)$ is also true.
Using the inequality of Proposition \ref{2.2} $\bf{(B)}$ for any integer $n$ we have
$$\lf\Vert\sum\li_{i=1}^{2^n}\lf|\al_i^{(n)}(m)\ri|\ph_i\ri\Vert
\ge{\frac{\sum\li_{i=1}^{2^n}\lf|\al_i^{(n)}(m)\ri|}{2^{n+1}}}\,
\lambda(2^n) \quad {\text{for any}} \quad m=1,2,\ld,2^n,$$
where the numbers $\al_i^{(n)}(m)$ are defined by (3.4). Hence for any integer $n$ we get
$$\max\li_{1\le m\le2^n}\lf\Vert\sum\li_{i=1}^{2^n}\lf|\al_i^{(n)}(m)\ri|\ph_i
\ri\Vert \ge{\frac{\max\li_{1\le
m\le2^n}{\sum\li_{i=1}^{2^n}\lf|\al_i^{(n)}(m)
\ri|}}{2^{n+1}}}\,\lambda(2^n). \eqno (3.12)$$ We know that
$\lf\Vert\sum\li_{i=1}^{2^n}\lf|\al_i^{(n)}(m)\ri|\ph_i\ri\Vert=
\ro^{(n)}(m)$ and $\sum\li_{i=1}^{2^n}\lf|\al_i^{(n)}(m)\ri|$
is the value of $\ro^{(n)}(m)$ in the space $l_1$
with respect to the natural basis. Therefore, by Theorem \ref{3.5} we have
$$\max\li_{1\le m\le2^n}\sum\li_{i=1}^{2^n}
\lf|\al_i^{(n)}(m)\ri|={\frac{3n+7}{9}}\,2^n+(-1)^n\,{\frac{2}{9}}
\quad {\text{for any}} \quad n=1,2,\ldots .$$
Putting these expressions into (3.12) we complete  the proof by elementary
calculations.
\qed
\vskip 0.2cm

\begin{remark} \label{3.8}
{\rm If a basis $(\ph_i)$ of a space $X$ is in addition
symmetric, then using the inequality of Proposition \ref{2.2} $\bf{(A)}$
we can prove by analogy with Theorem \ref{3.6} that
$$\ro^{(n)}\ge\max\left\{{\frac{n+2}{3}} \, \lambda(2^n), \, 2^n\right\},
\quad n=1,2,\ldots .$$}
\end{remark}
\vskip 0.2cm

It follows from Theorem \ref{3.6} that in spaces of type $p, \, p>1,$
for sufficiently large $n$ the lower estimate $2^n$ is more precise
than ${\frac{n+2}{6}} \,\lambda(2^n)$,
because in such spaces we have $\lambda(2^n)\le T_p(X) \,2^{n/p}$.
Hence, the lower estimate ${\frac{n+2}{6}} \,\la(2^n)$ can compete with
$2^n$ in spaces of type 1.

The following example shows that beside $l_1$ there exist
Banach spaces different from $l_1$ with $\sup\li_n\ro^{(n)}/2^n=\infty$.
\vskip 0.2cm

\begin{example} \label{3.9}
{\rm Consider the real function
$f(t)={\frac{\sqrt{\log_25}}{5}}\,{\frac{t+4}
{\sqrt{\log_2(t+4)}}},\,t\ge1$. It is concave since for every $t \ge 1$
we have
$$f^{''}(t)={\frac{\sqrt{\log_25}}{10 \ln2}} \cdot {\frac{-\log_2(t+4)+3/
(2\ln2)}{(t+4)\log_2^{5/2}(t+4)}}\le0.$$
By Proposition \ref{2.1}$(c)$ the sequence $(\lambda_n)$ with $\lambda_n=f(n),\, n=1,2,\ldots,$
is concave. Therefore, there exists at least one Banach space $X$
having a symmetric basis $(\ph_i)$ with symmetric constant equal to $1$ such that
$\lambda(n)=\lf\Vert\sum\li_{i=1}^n\ph_i\ri\Vert=\lambda_n$ for every
$n=1,2,\ldots $ (see \cite{LT}, p. 120). Hence by Remark \ref{3.8} for any integer $n$ we have
$$\ro^{(n)}\ge{\frac{n+2}{3}}\cdot\lambda(2^n)={\frac{n+2}{3}}\cdot{\frac{\sqrt{\log_25}}{5}}\cdot
{\frac{2^n+4}{\sqrt{\log_2(2^n+4)}}}>$$
$$>{\frac{\sqrt{\log_25}}{15}}\cdot{\frac{n+2}{\sqrt{n+2}}}\cdot2^n \ge
{\frac{\sqrt{\log_25}}{15}}\cdot{\sqrt{n+2}}\cdot2^n.$$
The space $X$ is not isomorphic to $l_1$ since
$$\lim\li_{n\to\infty}\sup_n{\frac{\lambda(n)}{n}}=\lim\li_{n\to\infty}\sup_n\left({\frac{\sqrt{\log_25}}{5}}\cdot
{\frac{n+4}{n\,{\sqrt{\log_2(n+4)}}}}\right)=0.$$

In particular, it follows from the obtained estimate
that the type of $X$ does not exceed $1$ (cf. Corollary \ref{3.4}).}
\end{example}

\section{Hadamard matrices}

The main aim of this section is to clear up whether the estimates for 
Sylvester matrices found in Section 3 can be extended to general Hadamard matrices. 
Let $\mathcal H_n^{\scriptstyle all}$ be the set of all
Hadamard matrices of order $n,\,n\in\mathbb{N}_{\mathcal H}$. For a Hadamard matrix
${\mathcal H}_n=\left[h_{ki}^n\right]$ we consider the same numerical
characteristic $\varrho_{\mathcal
H_n}(m)=\left\Vert\sum\limits_{i=1}^n\left(\sum\limits_{k=1}^mh_{ki}^n\right)
\varphi_i\right\Vert,\, m=1,2,\ldots,n,$ where $(\varphi_i)$ is a
normalized basis of a Banach space $X$. Setting
$a_k=\sum\limits_{i=1}^nh_{ki}^n\varphi_i,$ we notice that
$$\varrho_{\mathcal
H_n}(m)=\left\Vert\sum\limits_{k=1}^ma_k\right\Vert. \eqno (4.1)$$
If $(\varphi_i)$ is an unconditional basis with
unconditional constant equal to $1$, then we have $\max\limits_{1\le m\le
n}\varrho_{{\mathcal H}_n}(m)\le\lambda(n) \, n\le n^2$ for any
${\mathcal H}_n\in{\mathcal H}_n^{\scriptstyle all}$.
\vskip 0.2cm

Finally we set $\varrho_{\mathcal H_n}=\max\limits_{1\le m\le
n}\varrho_{{\mathcal H}_n}(m)$ \ and \
$\varrho_n=\max\limits_{{\mathcal H}_n\in{\mathcal
H}_n^{\scriptstyle all}}\varrho_{{\mathcal H}_n}.$
\vskip 0.2cm

\begin{remark} \label{4.1}
{\rm Note that the characteristic $\varrho_{{\mathcal
H}_n}=\varrho({\mathcal H}_n)$ can be regarded as a norm of the
Hadamard matrix ${\mathcal H}_n$. Indeed, let us denote by $\mathbf{M}_n$
the vector space of all square matrices of order
$n,\, n\in \mathbb{N}_{\mathcal H}$, and let $X$ be a Banach space with a basis $(\ph_i)$.
One has $\mathcal H_n^{{\scriptstyle all}}\subset\mathbf{M}_n$.
Let $\mathcal T_n=\left[t_{ki}^{n}\right]\in\mathbf{M}_n$ be a matrix
and $\ro(\mathcal T_n)=\max\li_{1\le m\le n}\lf\Vert\sum\li_{i=1}^n
\lf(\sum\li_{k=1}^mt_{ki}^n\ri)\ph_i\ri\Vert.$ It is easy to see that
$\ro$ is a norm in $\mathbf{M}_n$ and with respect to this norm $\mathbf{M}_n$
is a Banach space.}
\end{remark}

The following theorem gives us the lower estimate for $\varrho_n$.
\vskip 0.2cm

\begin{theor} \label{4.2}
Let $X$ be a Banach space with a
normalized unconditional basis whose unconditional constant is
$1$. Then we have
$$\varrho_n\ge\max\left\{(1/\sqrt2)\,\lambda(n)\,\sqrt n,\,n\right\} \quad
{\text{for any}} \quad n\in\mathbb{N}_{\mathcal H}.$$
\end{theor}

\Pro \
If one of the columns of a Hadamard matrix $\mathcal H_n$ consists
of $+1$ only, then we have $\ro_{\mathcal H_n}(n)=n$ and
the inequality $\ro_n\ge n$ is evident.

Let $\mathcal H_n=\left[h_{ki}^n\right]$ be a Hadamard matrix
of order $n$ and $(r_k(t))_{k\le n}$ be a sequence of
Rademacher functions defined on the interval $[0,1]$. For every $t\in [0,1]$ the matrix
$\mathcal H_{n,t}=\left[h_{ki}^n\,r_k(t)\right]$
is also a Hadamard matrix such that $\ro_{\mathcal
H_{n,t}}=\max\li_{1\le m\le n}\ro_{\mathcal H_{n,t}}(m)=
\max\li_{1\le m\le n}\lf\Vert\sum\li_{k=1}^ma_k\,r_k(t)\ri\Vert,$
where $a_k=\sum\li_{i=1}^nh_{ki}^n\ph_i, \ k=1,2,\ld,n.$

Let
$\xi(t)=\sum\li_{i=1}^n\lf|\sum\li_{k=1}^n\lag\ph_i^*,a_k\,r_k(t)\rag\ri|\ph_i$.
Using the fact that $(\ph_i)$ is an unconditional basis with
unconditional constant equal to $1$, it is easy to see that
$$||\xi(t)||=\lf\Vert\sum\li_{i=1}^n\sum\li_{k=1}^n\lag\ph_i^*,a_k\rag\,r_k(t)
\ph_i\ri\Vert=$$
$$=\lf\Vert\sum\li_{k=1}^n\lf(\sum\li_{i=1}^n\lag\ph_i^*,a_k\rag\,
\ph_i\ri)\,r_k(t)\ri\Vert=
\lf\Vert\sum\li_{k=1}^na_k\,r_k(t)\ri\Vert$$ for every
$t\in[0,1]$. As the Rademacher functions are bounded, $||\xi(t)||$ is integrable 
with respect to the Lebesgue measure on
$[0,1]$. Hence, there exists the Pettis integral $\e\,\xi$ of the measurable
function $\xi$ and $\e\,||\xi||\ge ||\e\,\xi||$. It is easy to see that
$\e\,\xi=\sum\li_{i=1}^n\left(\e\,\lf|\sum\li_{k=1}^n\lag\ph_i^*,a_k\rag\,r_k(t)\ri|
\right)\ph_i$.

As the Rademacher functions are bounded, $\ro_{\mathcal H_{n,t}}$ is also
integrable with respect to the Lebesgue
measure on $[0,1]$, and using the Khintchine's inequality we have
$$\infty>\e\,\ro_{\mathcal H_{n,t}}=\e\,\max\li_{1\le m\le n}
\lf\Vert\sum\li_{k=1}^ma_k\,r_k(t)\ri\Vert\ge \e\,||\xi||\ge||\e\,\xi||\ge$$
$$\ge(1/{\sqrt2})\lf\Vert\sum\li_{i=1}^n\lf(\sum\li_{k=1}^n\lag\ph_i^*,a_k
\rag^2\ri)^{1/2}\ph_i\ri\Vert=(1/{\sqrt2})\lambda(n)\,\sqrt n,$$
where $(\ph_i^*)$ are the biorthogonal functionals associated to the basis
$(\ph_i)$. Then, clearly, there exists a point $t_0\in [0,1]$ such that
$\ro_{\mathcal H_{n,t_0}}\ge \e\,\ro_{\mathcal H_{n,t}}$ and therefore
$\ro_n\ge\ro_{\mathcal H_{n,t_0}}\ge(1/{\sqrt2})\lambda(n)\,\sqrt n.$
\qed
\vskip 0.2cm

An immediate consequence of this theorem is the following corollary.
\vskip 0.2cm

\begin{co} \label{4.3}
In $l_p, \, 1\le p<2,$ with the natural basis 
we have $\sup\li_{n\in\mathbb{N}_{\mathcal H}}\ro_n/n=\infty$.
\end{co}
\vskip 0.2cm

For the spaces $l_p$ the similar fact for the Sylvester matrices holds
only for the space $l_1$ (see Theorem \ref{3.5}).

Let us estimate $\ro_n$ from above for the case of $l_p, \ 1\le p<\infty$.
\vskip 0.2cm

\begin{theor} \label{4.4}
In $l_p, \, 1\le p<\infty,$ with the
natural basis for any $n\in\mathbb{N}_{\mathcal H}$ the following inequality holds
$$\ro_n\le\max\left\{n^{(p+2)/2p},\,n\right\}.$$
\end{theor}

\Pro \
Let $p\ge2$ and $\mathcal H_n\in\mathcal H_n^{{\scriptstyle all}}$ be an 
arbitrary Hadamard matrix of order $n$. Using definition (4.1) and
the fact that $||a||_{l_p}\le||a||_{l_2}$ one can see that
$$\ro_{\mathcal H_n}\le\max\li_{1\le m\le n}
\lf\Vert\sum\li_{k=1}^ma_k\ri\Vert_{l_2}=\max\li_{1\le m\le n}\lf(\sum\li_{k=1}^ma_k,
\sum\li_{k=1}^ma_k\ri)^{1/2}=n,$$
where $(\cdot,\cdot)$ denotes the inner product in the space $l_2$.
Hence in $l_p, \ p\ge2,$ the estimate $\ro_n\le n$ holds.

Now let $1\le p \le2$ and $\mathcal H_n\in\mathcal H_n^{{\scriptstyle all}}$ be 
again an arbitrary Hadamard matrix of order $n$. If $a=(\al_i)\in l_p$
is a sequence of the length $n$ (i.e. $\al_n\ne0$ and $\al_i=0$ for any $i>n$),
then we have $||a||_{l_p}\le n^{(2-p)/2p}\,||a||_{l_2}.$ Hence, we have
$$\ro_{\mathcal H_n}\le n^{(2-p)/2p} \max\li_{1\le m\le n}\lf\Vert\sum\li_{k=1}^ma_k\ri\Vert_{l_2}=n^{(p+2)/2p}$$
and the theorem is proved.
\qed
\vskip 0.2cm

For Sylvester matrices Corollary \ref{4.3} and Theorem \ref{4.4}
yield the following corollary.
\vskip 0.2cm

\begin{co} \label{4.5}
Let $S^{(n)}$ be the Sylvester matrix of order
$2^n,\,n=1,2,\ldots$. Then in $l_p,\,p\ge2$, with the natural basis we have
$$\ro^{(n)}=2^n.$$
\end{co}
\vskip 0.2cm

Theorem \ref{4.2} and \ref{4.4} imply the following assertion.
\vskip 0.2cm

\begin{co} \label{4.6}
In $l_p,\,1\le p\le \infty$, with respect to
the natural basis for every $n\in\mathbb{N}_{\mathcal H}$ we have
$$(1/\sqrt2)\,n^{(p+2)/2p}\le\ro_n\le n^{(p+2)/2p}, \quad {\text{\it for}}
\quad 1\le p<2,$$
$$\ro_n=n, \quad {\text{\it for}} \quad p\ge2.$$
\end{co}
\vskip 0.2cm

Let $X$ be a Banach space (not necessarily with a basis),
$x_1,x_2,\ldots,x_n$ be a sequence of elements from the unit ball
of $X$ and ${\mathcal H}_n\in\mathcal H_n^{\scriptstyle
all},\,n\in\mathbb{N}_{\mathcal H}$. Let us put
$\hat\varrho_{\mathcal H_n}(m)=\left\Vert\sum\limits_{i=1}^n\left(\sum\limits_{k=1}^mh_{ki}^n\right)
x_i\right\Vert,\,m=1,2,\ldots,n$,\, $\hat\varrho_{\mathcal
H_n}=\max\limits_{1\le m\le n}\hat\varrho_{{\mathcal H}_n}(m)$
and $\hat\varrho_n=\max\limits_{{\mathcal H}_n\in\mathcal
H_n^{\scriptstyle all}}\hat\varrho_{{\mathcal H}_n}$.
\vskip 0.2cm

\begin{co} \label{4.7}
For any $n\in\mathbb{N}_{\mathcal H}$ we have
$\hat\varrho_n\le n\,\sqrt n.$
\end{co}

\Pro \ Using Corollary \ref{4.6} for the case $p=1$ the proof goes analogously to the proof of Corollary \ref{3.3}.
\qed
\vskip 0.2cm

Now we prove the analogue of Theorem \ref{3.1} for the Hadamard matrices.
\vskip 0.2cm

\begin{theor} \label{4.8}
Let $X$ be a Banach space with a normalized subsymmetric basis
whose subsymmetric constant is 1. Then we have for any $n\in\mathbb{N}_{\mathcal H}$
$$\varrho_n\le\lambda([\sqrt n]+1)\,n,$$
where $[\sqrt n]$ is the integer part
of $\sqrt n$.
\end{theor}

\Pro \
Let $\mathcal H_n=\left[h_{ki}^n\right]$ be a Hadamard matrix of order $n$.
As we already have noted
$$\ro_{\mathcal H_n}=\max\li_{1\le m\le n}\lf\Vert\sum\li_{i=1}^n\left(\sum\li_{k=1}^m
h_{ki}^n\right)\ph_i\ri\Vert\le\max\li_{1\le m\le
n}\sum\li_{i=1}^n\lf|\sum \li_{k=1}^mh_{ki}^n\ri|\le n\,\sqrt n
\eqno (4.2)$$ for every $\mathcal H_n\in\mathcal
H_n^{{\scriptstyle all}}$. For the sake of convenience let us introduce the notation
$$\al_i^{(n)}(m)=\lf|\sum\li_{k=1}^mh_{ki}^n\ri| \quad {\textit{for any}} \quad
i,m=1,2,\ld,n.\eqno (4.3)$$ Using the definition of the Hadamard matrices and (4.2)
we obtain the following properties of the numbers $\al_i^{(n)}(m)$:

$(a).$ For all $i$ and $m$ the number $\al_i^{(n)}(m)$ is an integer and
$0\le\al_i^{(n)}(m)\le n.$

$(b).$ For any $m$ we have $\sum\li_{i=1}^n\al_i^{(n)}(m)
\le n\,\sqrt n$.

Denote by $M$ the subset of $X$ consisting of $n$ points
$\{\sum\li_{i=1}^n\al_i^{(n)}(m)\ph_i: \, m=1,2,\ldots,n\},$ where $\al_i^{(n)}(m)$
is defined by (4.3). Then we have $\ro_{\mathcal H_n}=\max\li_{x\in M}||x||.$

Let us consider the following subsets of $X$:
$$S=\{\sum\li_{i=1}^nt_i\ph_i: 0\le t_i\le n, i=1,2,\ldots,n\} \quad
{\hbox{and}} \quad T=\{\sum\li_{i=1}^n t_i\ph_i:
\sum\li_{i=1}^n t_i\le n \sqrt n\}.$$
Since $S$ is an $n$-dimensional parallelepiped and $T$ is a hyperplane
in $X$, the sets $S,\,T$ as well as their intersection $S\cap T$ are convex.
Moreover, we have $M\subset S\cap T$. The set $S\cap T$ is compact because it is a bounded set 
in an $n$-dimensional subset of $X$ spanned by the basis vectors $\ph_1,\ph_2,\ld,\ph_n$. 
According to the Krein-Milman theorem
(see, for example, \cite{D}, p. 104) $S\cap T$ is a closed convex span of its extreme points. Hence we have
$$\ro_{\mathcal H_n}=\max\li_{x\in M}||x||\le\sup\li_{x\in S\cap T}||x||=
\sup\li_{x\in E}||x||, \eqno(4.4)$$
where $E$ is the set of all extreme points of $S\cap T$. The extreme points of
the set $S$ are the vertices of the parallelepiped $S$, i.e. the points of
the form $\sum\li_{i=1}^n\be_i\ph_i$, where each $\be_i$ takes the values
$0$ or $n$. Since $E \subset S\cap T$, the set $E$ contains those extreme 
points of $S$ for which the condition
$\sum\li_{i=1}^n\be_i\le n\,\sqrt n$ is satisfied. If we denote by $l$ the
number of these $\be_i$-s which are different from zero, then the last condition
can be expressed as follows: $ln\le n\,\sqrt n,$ or equivalently $l\le\sqrt n$.
Since $l$ is an integer, we get $l\le[\sqrt n].$ Since the
basis $(\ph_i)$ is subsymmetric, the norm $\lf\Vert\sum\li_{i=1}^n\be_i\ph_i\ri\Vert$ can be
estimated as follows
$\lf\Vert\sum\li_{i=1}^n\be_i\ph_i\ri\Vert\le\lambda(l)\,n<\lambda([\sqrt n]+1)\,n.$

It is easy to check that the set $E$, besides the vertices of the parallelepiped $S$,
contains the points of the intersection of the bound of $T$ with the edges of the
parallelepiped $S$. The edges of $S$ consists of the points which have the form $\sum\li_{i=1}^n\be_i\ph_i$, where
one of $\be_i$ satisfies the condition $0\le\be_{i_0}\le n$ and  all other $\be_i$-s take the
values $0$ or $n$. Denote by $l$ the number of $\be_i$-s for which $\be_i=n$. 
Due to the condition $\sum\li_{i=1}^n\be_i\ph_i\in T$,
we have $\be_{i_0}+ln\le n\,\sqrt n$. As $\be_{i_0}\ge 0$ and $l$
is an integer we have $l\le [\sqrt n]$. Since
$0\le\be_{i_0}\le n$, using again that $(\ph_i)$ is a subsymmetric basis, we obtain
$$\lf\Vert\sum\li_{i=1}^n\be_i\ph_i\ri\Vert=\lf\Vert\be_{i_0}\ph_{i_0}+\sum
\li_{i_0\ne i=1}^n\be_i\ph_i\ri\Vert\le\lf\Vert n\ph_{i_0}+\sum
\li_{i_0\ne i=1}^n\be_i\ph_i\ri\Vert\le\lambda([\sqrt n]+1)\,n.$$

Thus, for every point $x$ of the set $E$ the estimate $||x||\le\lambda([\sqrt n]+1)\,n$
is valid and using (4.4) we complete the proof of the theorem.
\qed
\vskip 0.2cm

\begin{remark} \label{4.9}
{\rm We can rephrase Theorem \ref{4.8} in the following way:
Let ${\mathcal H}_n=\left[h_{ki}^n\right]$ be a Hadamard matrix
of order $n\in\mathbb{N}_{\mathcal H}$ and let
$a_k=\sum\limits_{i=1}^nh_{ki}^n\varphi_i$,\,$k=1,2,\ldots,n,$
where $(\varphi_i)$ is a normalized subsymmetric basis of a Banach
space $X$ with subsymmetric constant equal to $1$. Then we have
$$\max\limits_{1\le m\le n}\left\Vert\sum\limits_{k=1}^m\vartheta_ka_k\right\Vert\le\lambda\left([{\sqrt n}]+1\right)\,n$$ for every sign
$\vartheta_k \in \{-1, 1 \},\,k=1,2,\ldots,n,$ every Hadamard matrix
${\mathcal H}_n\in{\mathcal H}_n^{\scriptstyle{all}}$ and every positive
integer $n\in\mathbb{N}_{\mathcal H}$.}
\end{remark}
\vskip 0.2cm

By Theorem \ref{3.1}, in a Banach space with a normalized subsymmetric
basis whose subsymmetric constant is $1$ we have
$\ro^{(n)}/\left(n\cdot2^n\right)\le 1.$ On the other hand, by Theorem \ref{3.5} in the
space $l_1$ we have $\ro^{(n)}/\left(n\cdot2^n\right)\ge1/3.$ Using Sylvester and
Hadamard matrices we can characterize the spaces isomorphic to $l_1$ as follows.
\vskip 0.2cm

\begin{theor} \label{4.10}
Let $X$ be a Banach space with a normalized subsymmetric basis $(\varphi_i)$ whose
subsymmetric constant is 1. The following statements are equivalent:

$(i)$. There is a constant $\delta>0$ such that
$\varrho_n/\left(n\,\sqrt n\right)\ge\delta$ for every
$n\in\mathbb{N}_{\mathcal H}$, where $\delta$ is independent of $n$.

$(ii).$ $X$ is isomorphic to $l_1$.

$(iii).$ There exists a constant $\eps>0$ which does not depend on $n$
such that for every $n=1,2,\ldots$ we have
$\ro^{(n)}/\left(n\cdot2^n\right)\ge\eps$.
\end{theor}

\Pro \
$(i)\Rightarrow(ii)$. Using Theorem \ref{4.8}
for every $n\in\mathbb{N}_{\mathcal H}$ we have
$$0<\delta\le\ro_n/(n\,\sqrt n)\le\lambda([\sqrt n]+1)\,n/(n\,\sqrt n)=
\lambda([\sqrt n]+1)/\sqrt n.$$
Therefore one has $\lambda([\sqrt n])/\sqrt n\ge\delta/2>0$ for infinitely many $n.$
Now the validity of the statement $(ii)$ follows from the fact
which was mentioned in Section 2: if
$$\lim\li_{n\to\infty}\sup\lambda(n)/n>0,$$
then $X$ is isomorphic to $l_1$.

$(ii)\Rightarrow(iii)$. Let $X$ be
isomorphic to $l_1$, and denote by $T:X\rightarrow l_1$ an isomorphism between
$X$ and $l_1$. It is clear that
$(T\ph_i)$ is an unconditional basis in $l_1$. Since in $l_1$ all normalized
unconditional bases are equivalent (see \cite{LT}, p. 71), there exists a
bounded linear operator $S:l_1\rightarrow l_1$ with bounded inverse operator,
such that $T\ph_i=Se_i$ for every integer $i$, where $(e_i)$ is a sequence of the
unit vectors in $l_1$. By Theorem \ref{3.5} for every integer $n$ we have
$$1/3\le\max\li_{1\le m\le2^n}\lf\Vert\sum\li_{i=1}^{2^n}\lf|\sum_{k=1}^m
s_{ki}^{(n)}\ri|e_i\ri\Vert/\left(n\cdot2^n\right)=$$
$$=\max\li_{1\le m\le2^n}\lf\Vert
\sum\li_{i=1}^{2^n}\lf|\sum\li_{k=1}^m s_{ki}^{(n)}\ri|S^{-1}T\ph_i\ri\Vert/
\left(n\cdot2^n\right)\le\Vert S^{-1}T\Vert
\max\li_{1\le m\le2^n}\ro^{(n)}(m)/\left(n\cdot2^n\right).$$
With $\eps=1/\left(3||S^{-1}T||\right)>0$ we get the validity
of assertion $(iii)$.

The implication $(iii)\Rightarrow (i)$ is true because $2^n\in\mathbb{N}_{\mathcal H}$.
\qed
\vskip 0.2cm

\section{Unsolved problem}

Let $(e_i)$ be the natural basis of the space $l_1$,
$\mathcal S^{(n)}=\left[s_{ki}^{(n)}\right]$ be the Sylvester matrix
of order $2^n$, $n=1,2, \ldots ,$ and $(a_k)_{k\le2^n}$ be the sequence in $l_1$ defined by
$$a_k=\sum\li_{i=1}^{2^n} s_{ki}^{(n)}e_i, \quad k=1,2,\ldots,2^n.$$
Let us formulate the assertion of Theorem \ref{3.5} in the following manner:
$$\varrho^{(n)}=\left\Vert\sum\limits_{k=1}^{m_n} a_k
\right\Vert_{l_1}={(3n+7)2^n}/9+{2(-1)^n}/9,$$ where $m_n={(2^{n+1}+(-1)^n)}/3$.

Now let us consider a permutation $\sigma:\{1,2,\ldots,2^n\}\to\{1,2,\ldots,2^n\}$ and
the following expression:
$$\left\Vert\sum\limits_{k=1}^{m_n} a_{\sigma(k)}\right\Vert_{l_1}.$$

By Corollary \ref{4.6} for every permutation $\sigma:\{1,2,\ldots,2^n\}\to\{1,2,\ldots,2^n\}$
we have
$$\left\Vert\sum\limits_{k=1}^{m_n} a_{\sigma(k)}\right\Vert_{l_1}\le 2^{3n/2}.$$

The authors do not know yet the answer for the following conjecture:
\vskip 0.2cm

\begin{conj} \label{5.1}
For any positive integer $n$ and for any permutation
$\sigma:\{1,2,\ldots,2^n\}\to\{1,2,\ldots,2^n\}$ the following inequality holds:
$$\left\Vert\sum\limits_{k=1}^{m_n} a_{\sigma(k)}\right\Vert_{l_1}
\ge{(3n+7)2^n}/9+{2(-1)^n}/9.$$
\end{conj}

{\sc\bf Acknowledgment.} The authors are sincerely grateful to the referees
for useful remarks that clearly improved the reading of this paper.

Authors addresses: \'A. Figula\\
    Department of Mathematics\\
    University of Debrecen\\
    H-4010 Debrecen, P.O. Box 12\\
    Hungary, figula@science.unideb.hu

\medskip
\noindent 
V. Kvaratskhelia\\
    Muskhelishvili Institute of Computational Mathematics\\
    of the Georgian Technical University,\\
    Akuri str., 8, Tbilisi 0160\\
    Sokhumi State University, A. Politkovskaia str., 9,\\
    Tbilisi 0186, Georgia, v\_kvaratskhelia@yahoo.com

\end{document}